\documentclass[12pt]{article}
\usepackage[koi8-r]{inputenc}
\usepackage[english]{babel}
\usepackage{amssymb}
\usepackage{amsmath}
\begin{document}
\title 
{On sharp constants in one-dimensional\\ embedding theorems of arbitrary order\footnote{Translation of the paper published in Russian in the collected volume 
``Problems of contemporary approximation theory'', St.Petersburg Univ. Publishers, 2004, 146--158. Some comments are added. I am greatful to Sergey Kryzhevich
for some bits of advice which helped to improve English text.}
}

\author{Alexander I. Nazarov\footnote{St.Petersburg Department of Steklov Institute, Fontanka, 27, St.Petersburg, 191023, Russia and 
St.Petersburg State University, Universitetskii pr. 28, St.Petersburg, 198504, Russia. 
Supported by St.Petersburg University grant 6.38.670.2013.
}
}

\date{}

\maketitle

\hfill {\it Dedicated to the memory}

\hfill {\it of Professor G.I.~Natanson}

\section{Statement of the problem and main results}

Let $r>0$ and $T>0$. Denote by $H^r(-T,T)$ the space of real $2T$-periodic functions with finite norm
$$\Vert y\Vert_{H^r}^2 =\sum_{k\in\mathbb Z}|\widehat y_k|^2\cdot\Big[\Big(\frac{\pi k}{T}\Big)^{2r}+1\Big],
$$
where $\{\widehat y_k\}$ are Fourier coefficients of the function $y$ with respect to exponential system on $(-T,T)$.\medskip

We consider the problem of finding the sharp (exact) constant in the embedding theorem $H^r(-T,T)\to L_q(-T,T)$, $2<q <\infty$:
$$\min\limits_{H^r}\frac {\Vert y\Vert_{H^r}}{\Vert y\Vert_{L_q}}=: \lambda_{q,r}(T) > 0.
\eqno(1.1)
$$

{\bf Remark\ 1}. If $q\le 2$ then due to the evident estimate $\Vert y\Vert_{L_q}\le(2T)^{1/2-1/q}\cdot\Vert y\Vert_{L_2}$ the function $c(x)\equiv 1$
provides the minimum in (1.1), and $\lambda_{q,r}(T)= (2T)^{1/2-1/q}$.\medskip

{\bf THEOREM\ 1}. Let $q>\left(\frac{\pi}{T}\right)^{2r}+2$. Then the inequality
$\lambda_{q,r}(T) <(2T)^{1/2-1/q}$ holds, and thus the function $c(x)\equiv1$ does not provide the minimum in (1.1).\medskip

{\it Proof}. Similarly to \cite{N}, we consider the following functional on $H^r(-T,T)$:
$$J_{q,r}(y) = \Vert y\Vert_{H^r}^2 -(2T)^{1-2/q}\Bigl ( \int\limits_{-T}^T 
\vert y\vert ^q\ dx\Bigr)^{2/q}.
$$
This functional evidently vanishes on the function $c$. Next, direct calculation shows that the first derivative
$DJ_{q,r}(c\,; h)$ also vanishes for arbitrary variation $h\in H^r(-T,T)$ while the second derivative is as follows:
$$\frac 12 D^2J_{q,r}(c; h,h) = \Vert h_1\Vert_{H^r}^2-\int\limits_{-T}^T (q-1) h_1^2\ dx,
\eqno(1.2)
$$
where
$$h_1(x) = h(x) - \frac{1}{2T}\int\limits_{-T}^T h\ dx.
$$

It is easy to see that (1.2) is just the quadratic form of the operator
$$ \left(- \frac {d^2}{dx^2}\right)^r - (q-2)
$$
at the interval $(-T,T)$ under periodic boundary conditions. Since $\int_{-T}^Th_1\ dx =0$ by construction, the smallest eigenvalue of this
operator is equal to $\left(\frac{\pi}{T}\right)^{2r}-(q-2)$. Therefore, for $q>\left(\frac{\pi}{T}\right)^{2r}+2$ the function $c$ does not provide
even local minimum to the functional $J_{q,r}$, and the statement follows.\hfill$\square$\medskip

{\bf Remark\ 2}. In the case $q=\infty$, $r\in\mathbb N$ sharp constants in (1.1) were calculated explicitely in \cite{GO}.\medskip

{\bf THEOREM\ 2}. Let $r\ge 1$ and $q\le\left(\frac{\pi}{T}\right)^{2r}+2$. Then the function $c(x)\equiv 1$ provides the minimum in (1.1), and thus 
$\lambda_{q,r}(T) =(2T)^{1/2-1/q}$.\medskip

{\bf Remark\ 3}. In \cite{N} Theorem 2 was proved for $r=1$ and $T=1$. However, the argument of \cite{N} cannot be ``scaled'' for the case of arbitrary $T$.\medskip

{\bf Remark\ 4}. In particular, this result refutes Theorem 2 in \cite{K}.\medskip

{\it Proof}. {\bf 1}. First, we assume $r=1$. In this case (1.1) can be rewritten as follows:
$$\lambda^2_{q,1}(T)=\min\limits_{H^1}\frac{\int\limits_{-T}^T(
y^{\prime 2}+y^2)\ dx} {\bigl ( \int\limits_{-T}^T \vert y\vert ^q\ dx\bigr)^{2/q}}. 
\eqno (1.3) 
$$

Standard argument shows that the extremal function in (1.3) does not change the sign on $[-T,T]$ (to be definite we assume that $y>0$) 
and satisfies the Euler--Lagrange equation (here $\mu$ is the Lagrange multiplier)
$$-y''+y=\frac q2\mu y^{q-1}.
\eqno(1.4)
$$
Integrating (1.4) we obtain
$$y'^2 = y^2 - \mu y^q - c_1.
\eqno(1.5)
$$
Moreover, if $y$ is not a constant then the right-hand side of (1.5) should have two zeros corresponding to maximal and minimal values of
$y$ at the period. This implies $c_1>0$.

Thus, any non-constant periodic positive solution of equation (1.4) corresponds to the motion along an oval given by equation (1.5) in the phase
plane. Since the period of this motion should be equal to $2T$, we can write
$$\int\limits_{y_0}^{y_1}\frac{dy}{\sqrt{y^2-c_1-\mu y^q}}=T
$$
(here $y_0$ and $y_1$ are the roots of denominator).

The change of variable $y=t\sqrt{c_1}$ gives
$$I_q(\alpha):=\int\limits_{x_0}^{x_1}\frac{dt}{\sqrt{t^2-1-\alpha t^q}}=T, 
\eqno(1.6)
$$
where $x_0(\alpha) < x_1(\alpha)$ are the roots of denominator in (1.6) and $\alpha>0$.

It is easy to see that the function $I_q(\alpha)$ is defined on the interval $]0,\alpha^*(q)[$
where
$$\alpha^*(q)= \frac 2{q - 2}\left (\frac {q - 2}q\right )^{q/2}
$$
is given by condition $x_0=x_1$.\medskip

The statement of Theorem follows from Lemma which will be proved in Section 2.\medskip

{\bf THE MAIN LEMMA.} For $\alpha \in\ ]0,\alpha^*(q)[$ the following inequality holds:
$$dI_q/d\alpha < 0.
\eqno(1.7)
$$
Furthermore,
$$\lim\limits_{\alpha \to\alpha^*(q)}I_q(\alpha)=\frac{\pi}{\sqrt{q-2}}.
\eqno(1.8)
$$

{\bf Remark\ 5}. The inequality (1.7) was conjectured in \cite{N}.\medskip

It follows from relations (1.7) and (1.8) that if $q\le\left(\frac{\pi}{T}\right)^{2}+2$ then
$I_q(\alpha)>T$ for all $\alpha \in \  ]0, \alpha^*(q)[$. Therefore equation (1.6) has no solutions if
$q\le\left(\frac{\pi}{T}\right)^{2}+2$. Thus constant function is a unique $2T$-periodic solution of equation (1.4),
and Theorem 2 is proved for $r=1$.\medskip

{\bf 2}. Now we consider the case $r>1$. Let $y\in H^r(-T,T)$. We introduce the notation
$T_1:=T\left(\frac{T}{\pi}\right)^{r-1}$,
$u(x):=y(\frac{T}{T_1}x)$. By the Steklov--Wirtinger inequality,
$$\Vert y\Vert^2_{H^r(-T,T)}\ge \int\limits_{-T}^T \Big(\Big
(\frac{\pi}{T}\Big)^{2(r-1)} y'^2 + y^2\Big)\ dx=
\Big (\frac{\pi}{T}\Big)^{r-1}\Vert u\Vert^2_{H^1(-T_1,T_1)}.
$$

If $q\le\bigl(\frac{\pi}{T_1}\bigr)^{2}+2=\left(\frac{\pi}{T}\right)^{2r}+2$ then the first part of the proof implies
\begin{equation*}
\left (\frac{\pi}{T}\right)^{r-1}\Vert
u\Vert^2_{H^1(-T_1,T_1)}\ge \left (\frac{\pi}{T}\right)^{r-1}(2T_1)^{1-2/q}\, 
\cdot\Vert u\Vert^2_{L_q(-T_1,T_1)}=(2T)^{1-2/q}\Vert y\Vert^2_{L_q(-T,T)},
\end{equation*}
and the proof is complete.\hfill$\square$\medskip

Note that for $r<1/2$ the statement of Theorem 2 is not true since for sufficiently
large $q$ the space $H^r$ is not embedded into $L_q$ (thus, the minimum in (1.1) equals zero for any $T$).\medskip

{\it CONJECTURE}. The statement of Theorem 2 holds for $r\ge 1/2$.\medskip

The proof of Theorem 2 implies also the following statement.\medskip

{\bf THEOREM\ 3}. Let $k\in {\mathbb N}$. If 
$$\left(\frac{(k-1)\pi}{T}\right)^2+2<q\le\left(\frac{k\pi}{T}\right)^2+2
\eqno(1.9)
$$ 
then the equation
$$-y''+y=y^{q-1}
\eqno(1.10)
$$
has exactly $k$ non-equivalent $2T$-periodic positive solutions.\medskip

{\it Proof}. For any $q>2$ equation (1.10) has a unique constant positive solution
$y(x)\equiv 1$. Next, any non-constant periodic positive solution corresponds to the motion along an 
oval given by equation (1.5) in the phase plane (with $\mu=2/q$). In this process one revolution should
take time $2T/n$, $n\in {\mathbb N}$. Hence we derive similarly to the proof of Theorem 2
that the equality 
$$I_q(\alpha)=T/n
$$ 
holds for some $\alpha \in \ ]0,\alpha^*(q)[$.

It follows from (1.7), (1.8) and evident  relation $\lim\limits_{\alpha\to0}I_q(\alpha) =+\infty$ 
that for $q$ meeting the inequalities (1.9) this equation is solvable if $n\le k-1$. This gives $k-1$ 
non-equivalent non-constant positive solutions. Theorem is proved.\hfill$\square$

\section{Proof of the Main Lemma}

Some of Lemmas in this Section were proved in \cite{N}. For the reader's convenience
we give them with full proofs.\medskip

First, we observe that the change of variable $t=(u^2\alpha)^{-\frac{1}{q-2}}$
in the integral $I_q(\alpha)$ gives the identity
$$I_q(\alpha)=\frac{2}{q-2}\cdot I_{\frac{2q}{q-2}}(\alpha^{2/(q-2)})
$$
for all $ \alpha \in \ ]0, \alpha^*(q)[$. Therefore it is sufficient to prove monotonicity of
$I_q$ (the inequality (1.7)) only for $2<q\le 4$.

Let us denote
$$f(t):=t^2-1-\alpha t^q,\qquad f(x_0)=f(x_1)=0.
\eqno(2.1)
$$
Then
$$f'(t)=2t-\alpha q t^{q-1},\qquad f'(x_0)>0,\ f'(x_1)<0.
$$

Denote by $\widehat x$ (unique positive) point where $f'$ changes sign. Then we have
$$f''(t)=2-\alpha q(q-1) t^{q-2}, \qquad f''(\widehat x) = -2(q-2);
\eqno(2.2)
$$
$$f'''(t)=-\alpha q(q-1)(q-2) t^{q-3}<0.
\eqno(2.3)
$$

Also we write down two evident equalities:
$$tf'(t) = qf(t)+q-(q-2)t^2,
\eqno(2.4)
$$
$$tf'''(t) = (q-2)f''(t)-2(q-2).
\eqno(2.5)
$$

{\bf LEMMA 2.1.} For $ \alpha \in \ ]0, \alpha^*[$ the identity
$$\frac {dI_q}{d\alpha} = - 4q(q-1)\int\limits_ {x_0}^{x_1}\frac{\sqrt{f} f't^{q-3}\ dt}{\psi^2}, 
\eqno(2.6)
$$
holds with $\psi=f'^2-2ff''$.\medskip

{\it Proof}. We have
$$ I_q^{(\epsilon)} (\alpha) :=\int\limits_{x_0+\epsilon}^{x_1-\epsilon} \frac{dt}{\sqrt{f}}\ \stackrel
{\epsilon\to0} \longrightarrow \ I_q(\alpha),
$$
and convergence is uniform in any compact subset of the interval $]0,\alpha^*[$.

Furthermore,
$$\frac{dI_q^{(\epsilon)}}{d\alpha} = \frac{dx_1}{d\alpha}\cdot \left.\frac 1{\sqrt{f}}\,
\right\vert^{x_1-\epsilon} -\ \frac{dx_0}{d\alpha}\cdot\left.\frac 1{\sqrt{f}}\, \right\vert^{x_0+\epsilon} +\ 
\frac 12\int\limits_{x_0+\epsilon} ^{x_1-\epsilon} \frac {t^q\ dt}{f^{3/2}}.
$$

However, (2.1) implies
$$\frac{dx_k}{d\alpha}=\left.\frac{t^q}{f'}\,\right\vert^{x_k}=
\left.\frac{t^q f' - 2qt^{q-1}f}{f'^2-2ff''}\,\right\vert^{x_k},\quad k=0,1,
\eqno(2.7)
$$
and thus
$$\frac{dI_q^{(\epsilon)}}{d\alpha} = \left.\frac 1{\sqrt{f}}\ \frac{t^q f'-2qt^{q-1}f}{\psi}\,
\right\vert_{x_0+\epsilon}^{x_1-\epsilon}\ +\ O(\epsilon^{1/2})\ +\ 
\frac 12\int\limits_{x_0+\epsilon}^{x_1-\epsilon} \frac {t^q\ dt}{f^{3/2}}.
$$

Note that $\psi(x_0)=f'^2(x_0)>0$ and
$$\psi'= -2f\cdot f'''>0\qquad \mbox{in}\quad  ]x_0,x_1[.
$$
Hence $\psi >0$ in $ [x_0,x_1]$, and we can write
$$\frac{dI_q^{(\epsilon)}}{d\alpha}=\int\limits_{x_0+\epsilon}^{x_1-\epsilon}
\left[\frac d{dt}\left(\frac {t^q f' - 2qt^{q-1} f}{\sqrt{f} \psi}\right)
 +\ \frac{t^q }{2f^{3/2}}\right]\ dt \ +\ O(\epsilon^{1/2}).
 $$

The expression in square brackets is equal to
$$- 4q(q-1)\frac{\sqrt{f} f' t^{q-3}}{\psi^2}.
$$
Therefore $dI_q^{(\epsilon)}/d\alpha$ converges to the right-hand side of (2.6) as $\epsilon \to 0$. 
Moreover, convergence is uniform  in any compact subset of the interval $]0,\alpha^*[\,$. 
This completes the proof.\hfill$\square$\medskip

Let us denote $g(t):=\psi^2(t)/t^{q-3}$. Then formula (2.6) can be rewritten using the mean value theorem:
\begin{equation*}
\frac {dI_q}{d\alpha} = - 4q(q-1)\ \Bigl(\,\int\limits_{x_0}^{\widehat x}+\int\limits_{\widehat x}^{x_1}\,\Bigr)\,\frac{\sqrt{f} f'\ dt}g
= \frac {- 8q(q-1)}3\ f^{3/2}(\widehat x)\ \left( \frac1{g(t_0)} - \frac1{g(t_1)} \right),
\tag{2.8}
\end{equation*}
where $t_0$ and $t_1$ are some points in the intervals $]x_0,\widehat x[$ and $]\widehat x, x_1[\,$, respectively.

Since the function $\psi$ increases on the segment $[x_0,x_1]$, the function $g$ also increases on $ [x_0,x_1]$ if $q\le 3$. Thus,
$g(t_0)<g(t_1)$, and (2.8) implies (1.7) for $2<q\le 3$.\medskip

Now we turn to the case $3\le q\le4$.\medskip

{\bf LEMMA 2.2.} Let $q>3$. Then the function $g'$ changes the sign in $[x_0,x_1]$
exactly twice. Moreover,
$$g'(x_k)<0,\quad k=0,1;\qquad g'(\widehat x)>0.
\eqno(2.9)
$$

{\it Proof}. We observe that
$$g'(t) = \frac {\psi(t)}{t^{q-2}} \Bigl(2\psi'(t)t-(q-3)\psi(t)\Bigr).
$$

Let us denote by $g_1(t)$ the expression in large brackets. Then obviously
$$g_1(x_k)=-(q-3)\psi(x_k)<0,\qquad k=0,1;
$$ 
on the other hand, (2.2) and (2.5) give
$$g_1(\widehat x) = 4(q+1)(q-2) f(\widehat x) >0,
$$
and the inequalities (2.9) are proved. 

Further,
$$g_1'(t) = -2f'''(t) \Bigl( 2f'(t)t+(q-1)f(t)\Bigr).
$$
Let us denote by $g_2(t)$ the expression in large brackets. Then
$$\aligned
g_2'(t)&=(q+1)f'(t)+2tf''(t),\\
g_2''(t)&=(q+3)f''(t)+2tf'''(t)=(3q-1)f''(t)-4(q-2)\qquad\mbox{by (2.5)},\\
g_2'''(t) &= (3q-1)f'''(t) <0.
\endaligned
$$
Hence $g_2''$ changes sign at most once, and $g_2'$ changes sign at most
twice. Moreover, (2.2) and (2.3) imply $g_2''(x_1)<0$.

If $g_2''(x_0) >0$ then $f''(x_0)>0$ whence $g_2'(x_0)>0$. Since $g_2'(x_1)<0$ always, in this case $g_2'$ 
changes sign exactly once. Otherwise, if $g_2''(x_0) <0$ then $g_2''$ does not change sign and $g_2'$ changes sign at most once.

So, in any case $g_2$ changes sign at most twice. However, it is evident that $g_2(x_0)>0$ and $g_2(x_1)<0$. Thus $g_2$ changes sign exactly once,
and $g_1$ changes sign exactly twice. The Lemma is proved.\hfill$\square$\medskip

{\bf LEMMA 2.3.} Let $3\le q\le 4$. Then for $\alpha \in]0,\alpha^*(q)[$ the following inequalities hold:
$$g(x_0)<g(\widehat x) <g(x_1).
\eqno (2.10)
$$

{\it Proof}. First, we observe that for $\alpha=\alpha^*(q) $ the segment $[x_0,x_1]$ degenerates into the point
$$x^*=\sqrt{q/(q-2)}.
$$
Hence $x_0=x_1=\widehat x=x^*$, and $g(x^*)=0$.

Next, by (2.4) we have
$$g(x_k)=\beta(x_k),\quad k=0,1; \qquad g(\widehat x)= \gamma(\widehat x),$$
where
$$\beta(t) = \frac{\bigl((q-2)t^2-q\bigr)^4}{t^{q+1}}\,,\qquad
\gamma(t) = \left(\frac{16(q-2)}q\right)^2\,\frac {\bigl((q-2)t^2-q\bigr)^2} {t^{q-3}}.
$$
Hence
$$\frac d{d\alpha}\{g(x_k)\} = \beta'(x_k) \frac {dx_k}{d\alpha},\quad k=0,1; \qquad 
\frac d{d\alpha}\{g(\widehat x)\} = \gamma'(\widehat x)\frac{d\widehat x}{d\alpha}.
$$

Formulae (2.7) and (2.4) imply
$$\frac{dx_k}{d\alpha}=\left.\frac{t^{q+1}}{q-(q-2)t^2}\,\right\vert^{x_k};\qquad 
\frac{dx_0}{d\alpha}>0,\ \frac{dx_1}{d\alpha}<0.
\eqno(2.11)
$$
In a similar way, relations (2.2) and $f'(\widehat x)=0$ imply
$$\frac{d\widehat x}{d\alpha}=\left.\frac{-qt^{q-1}}{2(q-2)}\,\right\vert^{\widehat x}<0.
\eqno(2.12)
$$
Hence
$$\frac d{d\alpha}\{g(x_k)\} = \beta_1(x_k),\quad k=0,1; \qquad 
\frac d{d\alpha}\{g(\widehat x)\} = \gamma_1(\widehat x),
$$
where
\begin{equation*}
\aligned
\beta_1(t) &= \frac{\bigl((q-2)t^2-q\bigr)^2}t\Bigl((q-2)(q-7)t^2 -q(q+1)\Bigr),\vphantom{\frac{b^b}{g^g}} \\
\gamma_1(t) &= \frac{8(q-2)t}q\bigl((q-2)t^2-q\bigr)\Bigl((q-2)(q-7)t^2 -q(q-3)\Bigr),\vphantom{\frac{b^b}{g^g}}
\endaligned
\end{equation*}

In the same way we obtain for $j=2,3$
$$\frac {d^j}{d\alpha^j}\{g(x_k)\} = \beta_j(x_k),\quad k=0,1; \qquad 
\frac {d^j}{d\alpha^j}\{g(\widehat x)\} = \gamma_j(\widehat x),$$
where
\begin{equation*}
\aligned
\beta_2(t) = -t^{q-1}\bigl[5(q-7)(q-2)^2t^4-4q(q-1)(q-2)t^2-q^2(q+1)\bigr],
\\
\gamma_2(t) =  -4t^{q-1}\bigl[5(q-7)(q-2)^2t^4-6q(q-2)(q-5)t^2+q^2(q-3)\bigr],
\endaligned
\end{equation*}
\begin{multline*}
\beta_3(t) = \frac{t^{2q-1}}{(q-2)t^2-q}\\
\times\bigl[5(q+3)(q-7)(q-2)^2t^4
-4q(q-1)(q-2)(q+1)t^2-q^2(q-1)(q+1)\bigr],
\tag{2.13}
\end{multline*}
$$\gamma_3(t) = \frac {2t^{2q-3}}{q-2}\,\bigl[5(q+3)(q-7)(q-2)^2t^4
-6q(q-2)(q-5)(q+1)t^2+q^2(q-1)(q-3)\bigr].
\eqno(2.14)
$$ 

Note that 
$$\beta_1(x^*) = \gamma_1(x^*) = 0;\qquad \beta_2(x^*) = \gamma_2(x^*) = 32\,q^{2}\,(x^*)^{q - 1}.
\eqno(2.15)
$$
The inequalities (2.11) and (2.12) imply that $x_0<x^*<\widehat x$ for $0<\alpha<\alpha^*(q)$. Hence the relations (2.13) and (2.14) give for
$3\le q\le 4$
$$\beta_3(x_0)>0>\gamma_3(\widehat x).
$$
We integrate this inequality triply with respect to $\alpha$ using (2.15) and obtain the first inequality in (2.10).\medskip

To prove the second inequality in (2.10) we observe that for $q$ under consideration we have $\gamma_3'(t)<0$ for $t>x^*$ whence $\gamma_3(\widehat x)>\gamma_3(x_1)$.

Let us compare $\beta_3(x_1)$ and $\gamma_3(x_1)$. To proceed we introduce a new variable
$$z:=x_1\sqrt{(q-2)/q}
$$
and note that $z$ increases from $1$ to $+\infty$ as $\alpha$ decreases from $\alpha^*(q)$ to $0$. Furthermore, by direct calculation we obtain
$$\frac{\gamma_3(x_1)}{\beta_3(x_1)}-1=\frac{P(z,q)}{ z^{2}\bigl(5z^{4}(q^{2}-4q - 21)-(4z^{2} + 1)(q^{2}-1)\bigr)},
\eqno(2.16)
$$
where
\begin{equation*}
\aligned
P(z,q)=&-5(2q - 1)(7 - q)(q + 3)z^{6}-2(q + 2)(11q^{2} - 68q + 1)z^{4} \\
&+(14q^{3}-55q^{2}- 54q  - 1)z^{2} - 2q(q - 1)(q - 3).
\endaligned
\end{equation*}

It is easy to see that for $q$ under consideration the denominator in (2.16) is negative while the numerator
is positive for $z=1$ and decreases with respect to $z$. Furthermore,
numerical evaluation of zeros of the polynomial $P(1.15,q)$ gives
\begin{equation*}
P\,\bigr|_{z=1.15}=1.16747016\cdot(q + 5.796999289)\cdot(q - 4.076622243)\cdot(q -8.501415029).
\end{equation*}
This shows that the numerator in (2.16) is still positive for $z=z_0:=1.15$ and $3\le q\le 4$.
Thus, $\gamma_3(x_1)/\beta_3(x_1)<1$ for $z\in [1,z_0]$. Since $\beta_3(x_1)<0$, the inequality
$$\beta_3(x_1)<\gamma_3(x_1)<\gamma_3(\widehat x)
$$
holds for such $z$. We integrate this inequality triply with respect to $\alpha$ using (2.15) 
and obtain that the second inequality in (2.10) is fulfilled for $z\le z_0$.\medskip

To progress further, we introduce the notation
$$\tau:=\frac{\widehat x}{x_1}=\left(\frac{2z^{2}}{z^{2}q - (q - 2)}\right)^{\frac{1}{q - 2}}
$$
and note that $\tau$ increases with respect to $q$ and decreases with respect to $z$. 
Since $\tau\,\bigr|_{z=z_0,q=4}=0.8966333519$, we have $\tau<\tau_0:=0.897$ for $z>z_0$.

Furthermore, we observe that $\gamma_1'(t)<0 $ for $t>x^*$ and for $q$ under consideration,
whence $\gamma_1(\widehat x)>\gamma_1(\tau_0x_1)$ for $z>z_0$.

Let us compare $\beta_1(x_1)$ and $\gamma_1(\tau_0x_1)$. Direct calculation gives
$$\frac{\gamma_1(\tau x_1)}{\beta_1(x_1)}-1=\frac{Q(z,q,\tau)} { q(z -1)^{2}(z+ 1)^{2}(z^{2}(7-q) + (q + 1))},
\eqno(2.17)
$$
where
\begin{equation*}
\aligned
Q(z,q,\tau) =&-(7 - q)(q-8\tau ^{5})z^{6}+ ( q(13- 3q)  - 16\tau^{3}(5-q))z^{4} \\
&+(q(3q - 5) - 8\tau(q-3))z^{2}- q(q + 1).
\endaligned
\end{equation*}

The denominator in (2.17) is evidently positive for $3\le q\le4$. To estimate the
numerator we observe that $Q$ is convex with respect to $q$. Hence
$$\frac{\partial Q}{\partial q}\le\left.\frac{\partial Q}{\partial q}\,
\right|_{q=4}= -8\tau z^2(\tau^2z^2-1)^2-(9-z^2)(z^2-1)^2,
\eqno(2.18)
$$
and therefore the polynomial $Q$ decreases with respect to $q$ at least for $z\le 3$.

Numerical evaluation of zeros of the polynomial $Q(z,3,\tau_0)$ gives
\begin{equation*}
\aligned
Q\,\bigr|_{q=3,\tau=\tau_0} &=6.582844536\,\cdot(z + 1.157682736)\cdot(z -1.157682736)\\
&\times(z^{2} + 1.636399020z  + 1.166257009)\cdot(z^{2} - 1.636399020z +1.166257009).
\endaligned
\end{equation*}
This shows that the numerator in (2.17) is negative for $z\le z_1:=1.157$,
$\tau=\tau_0$ and $3\le q\le4$. Thus, $\gamma_1(\tau_0x_1)/\beta_1(x_1)<1$ for $z\in [z_0,z_1]$. 
Since $\beta_1(x_1)<0$, the inequality
$$\beta_1(x_1)<\gamma_1(\tau_0x_1)<\gamma_1(\widehat x)
$$
holds for such $z$. We integrate this inequality with respect to $\alpha$  
and obtain that the second inequality in (2.10) is fulfilled for $z\le z_1$.\medskip

The next step is completely similar. Since 
$\tau\,\bigr|_{z=z_1,q=4}=0.8933635819$, the inequality
$\tau<\tau_1:=0.894$ holds for $z>z_1$.

Numerical evaluation of zeros of the polynomial $Q(z,3,\tau_1)$ gives
\begin{equation*}
\aligned
Q\,\bigr|_{q=3,\tau=\tau_1} &=6.274166280\cdot(z+ 1.166008256)\cdot(z -1.166008256)\\
&\times(z^{2} + 1.656562646z  + 1.186071814)\cdot(z^{2} - 1.656562646z+1.186071814).
\endaligned
\end{equation*}
This shows that the numerator in (2.17) is negative for $z\le z_2:=1.166$,
$\tau=\tau_1$ and $3\le q\le4$. Hence the inequality
$$\beta_1(x_1)<\gamma_1(\tau_1x_1)<\gamma_1(\widehat x)
$$
holds for $z\in [z_1,z_2]$. Integrating this inequality with respect to $\alpha$
we obtain the second inequality in (2.10) for $z\le z_2$.

Further, by the same way we obtain for $z>z_k$ the estimate $\tau<\tau_k$ which helps
to stretch out the second inequality in (2.10) to $z=z_{k+1}$. Thus we obtain the chain
$$\begin{array} {ll}
z_2=1.166 &\tau_2=0.890\\ 
z_3=1.177 &\tau_3=0.885\\ 
z_4=1.194 &\tau_4=0.878\\ 
z_5=1.221 &\tau_5=0.868\\
z_6=1.271 &\tau_6=0.851\\
z_7=1.438 &\tau_7=0.815
\end{array}
$$

Finally, substituting $\tau=\tau_7$ into (2.18) one can make sure that $\partial Q/\partial q\,\bigr|_{\tau=\tau_7}<0$ 
for all $z$. Therefore,
\begin{equation*}
\aligned
Q\,\bigr|_{\tau=\tau_7} &\le Q\,\bigr|_{q=3,\tau=\tau_7}=- 0.493638296\cdot(z^{2} +12.72622203) \\
&\times(z^{2} +2.170075688z+1.376692111)\cdot(z^{2} - 2.170075688z+1.376692111)<0.
\endaligned
\end{equation*}
This gives the second inequality in (2.10) for $z>z_7$. The proof is complete.\hfill$\square$\medskip

It follows from Lemmas 2.2 and 2.3 that  $g(t_0)<g(t_1)$ in (2.8) for $3\le q\le4$. This completes the proof of the inequality (1.7) for $2<q\le4$.\medskip

To prove (1.8) we observe that by Rolle's Theorem for any $t \in\ ]x_0,x_1[$ there exists $\widetilde t(t) \in \ ]x_0,x_1[$ such that
$$f(t) = - \frac {f''(\widetilde t)}2\, (x_1-t)(t-x_0).
$$
Hence
$$I_q(\alpha) = \sqrt{\frac {-2}{f''(\widehat t)}}\,\int\limits_{x_0}^{x_1}\frac
{dt}{\sqrt{(x_1-t)(t-x_0)}} =\pi \sqrt{\frac {-2}{f''(\widehat t)}},
$$
where $\widehat t$ is some point in $]x_0,x_1[$.

Furthermore, $\widehat t \to x^*=\widehat x$ as $\alpha \to \alpha^*$, and (2.2) implies
$$f''(\widehat t) \to f''(\widehat x) = -2(q-2),
$$ 
whence (1.8) follows.\hfill$\square$

\end{document}